\documentclass[11pt]{article}
\usepackage{amsmath,amssymb,amsthm}
\usepackage{caption}
\usepackage{derivative}
\usepackage[shortlabels]{enumitem}
\usepackage{epigraph}
\usepackage{float}
\usepackage[margin=1in]{geometry}
\usepackage[colorlinks]{hyperref}
\usepackage{listings}
\usepackage{mathtools}
\usepackage[framemethod=TikZ]{mdframed}
\usepackage{microtype}
\usepackage[section]{placeins}
\usepackage[headsepline]{scrlayer-scrpage}
\usepackage{setspace}
\usepackage{thmtools}
\usepackage[usenames,dvipsnames,svgnames]{xcolor}

\providecommand{\clubg}[1]{\bgroup\color{green!40!black}[$#1\clubsuit$]\egroup}

\DeclareMathOperator{\Av}{Av}

\DeclareMathOperator{\tl}{tl}
\DeclareMathOperator{\LRmax}{LRmax}
\DeclareMathOperator{\del}{del}
\reversemarginpar

\mdfdefinestyle{mdbluebox}{roundcorner=10pt,innerbottommargin=9pt,
    linecolor=blue,backgroundcolor=TealBlue!5,}
\declaretheoremstyle[headfont=\sffamily\bfseries\color{MidnightBlue},
    mdframed={style=mdbluebox},]{thmbluebox}
\mdfdefinestyle{mdredbox}{frametitlefont=\bfseries,innerbottommargin=8pt,
    nobreak=true,backgroundcolor=Salmon!5,linecolor=RawSienna,}
\declaretheoremstyle[headfont=\bfseries\color{RawSienna},
    mdframed={style=mdredbox},headpunct={\\[3pt]},postheadspace=0pt,]{thmredbox}
\mdfdefinestyle{mdgreenbox}{linecolor=ForestGreen,backgroundcolor=ForestGreen!5,
    linewidth=2pt,rightline=false,leftline=true,topline=false,bottomline=false,}
\declaretheoremstyle[headfont=\bfseries\sffamily\color{ForestGreen!70!black},
    mdframed={style=mdgreenbox},headpunct={ --- },]{thmgreenbox}
\mdfdefinestyle{mdblackbox}{linecolor=black,backgroundcolor=RedViolet!5!gray!5,
    linewidth=3pt,rightline=false,leftline=true,topline=false,bottomline=false,}
\declaretheoremstyle[mdframed={style=mdblackbox}]{thmblackbox}

\declaretheorem[style=thmbluebox,name=Theorem,numbered=no]{theorem*}
\declaretheorem[style=thmbluebox,name=Lemma,numbered=no]{lemma*}

\declaretheorem[style=thmgreenbox,name=Claim,numbered=no]{claim*}

\declaretheorem[style=thmblackbox,name=Remark,numbered=no]{remark*}

\declaretheorem[style=thmgreenbox,name=Definition,numbered=no]{definition*}

\declaretheorem[style=thmblackbox,name=Example,numbered=no]{example*}

\newlist{walk}{enumerate}{3}
\setlist[walk]{label=\bfseries (\alph*)}
\newenvironment{subproof}[1][Subproof]{%
\begin{proof}[#1] }%
{\end{proof}}

\usepackage{asymptote}
\begin{asydef}
size(10cm); // set a reasonable default
usepackage("amsmath");
usepackage("amssymb");
settings.tex="pdflatex";
settings.outformat="pdf";
import geometry;
void filldraw(picture pic = currentpicture, conic g, pen fillpen=defaultpen, pen drawpen=defaultpen) { filldraw(pic, (path) g, fillpen, drawpen); }
void fill(picture pic = currentpicture, conic g, pen p=defaultpen) { filldraw(pic, (path) g, p); }
pair foot(pair P, pair A, pair B) { return foot(triangle(A,B,P).VC); }
pair centroid(pair A, pair B, pair C) { return (A+B+C)/3; }
\end{asydef}

\begin{document}
\title{A Characterization of the $2m-4$ Case of Highly Sorted Permutations}
\author{Kai Yi}
\date{\today}
\maketitle

\begin{abstract}
Let $s$ denote West's stack-sorting map. In 2020, Defant characterized and enumerated the set $s^{n-m}(S_n)$ for $n \geq 2m-3$. While $|s^{n-m}(S_n)| = B_m$ when $n \geq 2m-2$, where $B_m$ denotes the $m$th Bell number, there are additional permutations when $n = 2m-3$. In this paper, we explore the more complex $n = 2m-4$ case, with several forms of additional permutations. We characterize $s^{m-4}(S_{2m-4})$ and find that its size is \[B_m + \frac{m^2 + 7m - 28}{2}\] for $m \geq 5$. This answers Defant's question about the $2m-4$ case. Furthermore, we find some differences in the behavior of the $2m-5$ case compared to the $2m-3$ and $2m-4$ cases.
\end{abstract}

\section{Introduction}

Stack sorting was first introduced in 1968 by Knuth in \cite{Knuth}, who proved a permutation can be sorted into the identity permutation using the stack sort machine (push and pop operations) if and only if it avoids the pattern $231$. The formal definition of pattern avoidance is:

Let $p$ and $q$ be permutations. Let the entries of $q$ be $q_1q_2 \dots q_k$. If there exists a subsequence of entries $p_1 p_2 \dots p_k$ in $p$ such that $p_i < p_j$ if and only if $q_i < q_j$ for all indices $i, j$, we say $p$ contains the pattern $q$. Otherwise, we say $p$ avoids the pattern $q$.

For example, $13425$ avoids the pattern $321$ but contains the pattern $123$ because of the entries $1, 3, 5$.
\newline

In 1990, West \cite{West} introduced the stack-sorting map $s : S_n \rightarrow S_n$, which ``pushes" an entry into a stack if the entries in the stack would remain in increasing order from top to bottom and ``pops" an entry from the top of the stack if this would not happen. This is similar to the Tower of Hanoi, where you cannot place a large disk on top of a small one.

For example, here is how the permutation $15243$ is sorted using West's function:
\begin{figure}[H]
    \centering
    \includegraphics[width=0.6\linewidth]{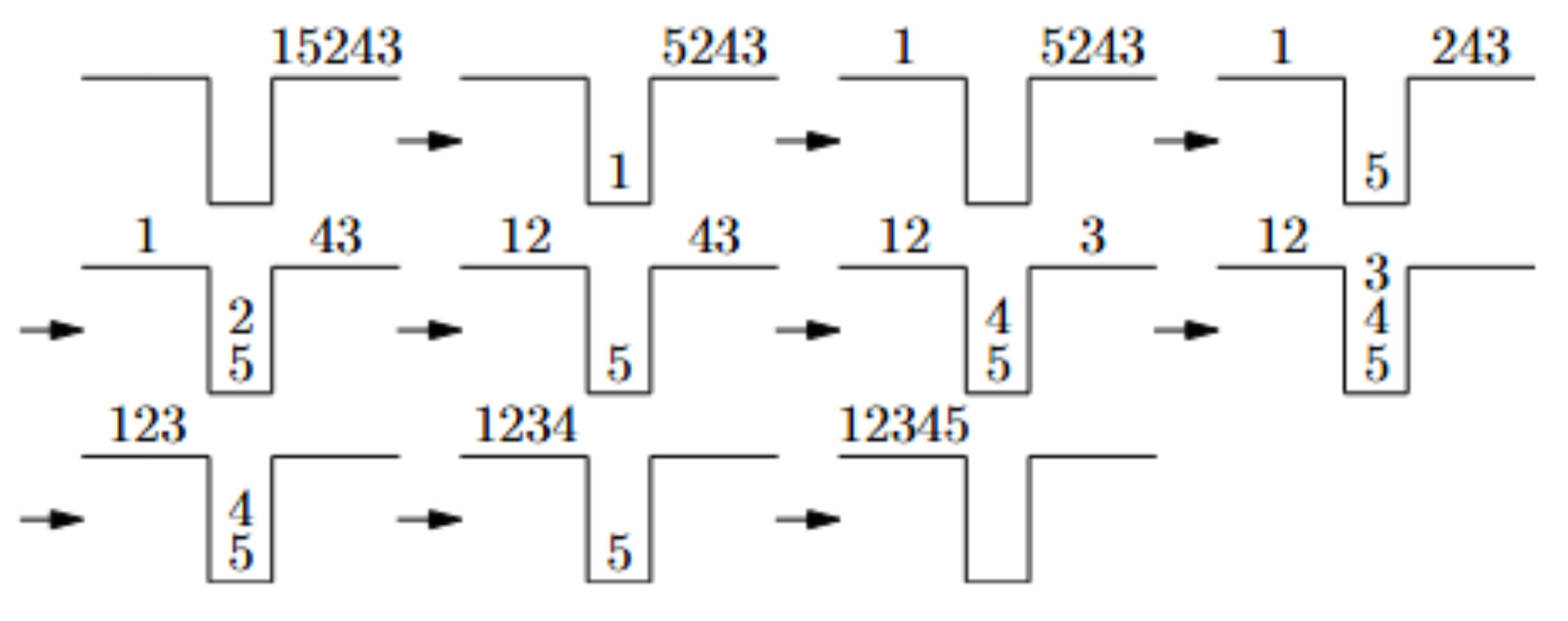}
    \caption*{\footnotesize Figure 1: Stack Sorting $15243$ using West's function}
\end{figure}

In the second step, we pop $1$ rather than push $5$ because pushing $5$ would result in $5$ (a large entry) being in the stack on top of $1$ (a small entry).
\newline 

While the effects of a single stack sort operation are simple to understand \cite{BonaWalk}, the effects of multiple stack sort iterations have been explored in several papers already \cite{Zeilberger} \cite{Bousquet} \cite{Claesson} \cite{DefantDescents}. West's stack sorting map attracts great interest from researchers because it is a rich combinatorial structure, with connections to trees, recursion, Catalan numbers, and more sophisticated integer sequences even though it was initially intended for sorting permutations.

One result regarding the effects of multiple stack sort operations is that the maximum number of stacks needed to sort a length $n$ permutation is $n-1$ \cite{BonaWalk}. For example, $23451$ is sorted using $n-1 = 4$ stacks.
\[23451 \xrightarrow{s} 23415 \xrightarrow{s} 23145 \xrightarrow{s} 21345 \xrightarrow{s} 12345\]

In 2020, Defant \cite{DefantBell} extends this result by studying the characterization and enumeration of highly-sorted permutations; in other words, the permutations in $s^t(S_n)$ when $t$ is close to $n$. His two main results are shown below. The relevant definitions follow in the next section.
\newline

\noindent
\textbf{Theorem 1.1:} Let $m$ and $n$ be positive integers such that $n \geq 2m-2$. A permutation $\pi \in S_n$ is $(n-m)$-sorted if and only if $\tl(\pi) \geq n-m$ and every descent top of $\pi$ is a left-to-right maximum of $\pi$. Consequently,
\[|s^{n-m}(S_n)| = B_m,\]
where $B_m$ denotes the $m$th Bell number.
\newline

\noindent
\textbf{Theorem 1.2:} Let $m \geq 3$ be an integer. A permutation $\pi \in S_{2m-3}$ is $(m-3)$-sorted if and only if one of the following holds:
\begin{itemize}
    \item $\tl(\pi) \geq m-3$ and every descent top of $\pi$ is a left-to-right maximum of $\pi$;
    \item $\pi = \zeta_{\ell, m}$ for some $\ell \in \{3, \dots, m\}$. The permutation $\zeta_{\ell, m}$ is obtained from the identity permutation of length $2m-3$ by swapping entries $1$ and $2$ and moving entry $\ell$ to the front.
\end{itemize}
Consequently,
\[|s^{m-3}(S_{2m-3})| = B_m + m-2.\]

When $m = n$, we have $|s^{n-m}(S_n)| = |s^0(S_m)| = |S_m| = m!$. Therefore, a natural question that arises is: when $n$ decreases past $2m-2$ and toward $m$, how does $|s^{n-m}(S_n)|$ increase from $B_m$ to $m!$? A natural next step is to first explore the $n = 2m-4$ case.

The main result we will prove is the following.
\newline

\noindent
\textbf{Theorem 1.3:} Let $m \geq 4$ be an integer. A permutation $\pi \in S_{2m-4}$ is $(m-4)$-sorted if and only if one of the following holds:
\begin{itemize}
    \item $\tl(\pi) \geq m-4$ and every descent top of $\pi$ is a left-to-right maximum of $\pi$.
    \item The first four (three for the last case) entries of $\pi$ are of the form $c312, 2c31, c231, 1c32, c132, dc21,$ or $d21$, where $4 \leq c \leq m$ and $3 \leq d \leq m$ and $d < c$ (each of the three inequalities is used if applicable), and the remaining entries are placed in increasing order.

    In all these cases, the first four entries of $\pi$ have exactly two descents, and neither $c$ nor $d$ is the entry for the second descent.
    \item $\pi = 4321$.
\end{itemize}
Consequently,
\[|s^{m-4}(S_{2m-4})| = B_m + \frac{m^2 + 7m - 28}{2}\]
if $m \geq 5$. Add one to the above value if $m = 4$.

The $2m-4$ case is less stable compared to the $2m-3$ case. Instead of having one exceptional form other than the Bell number case, there are now several exceptional forms. Furthermore, the number of permutations with these exceptional forms is quadratic instead of linear. After proving the $2m-4$ case, we also provide insights about the complexity of the $2m-5$ case.

\section{Preliminaries}

To begin, we will first present definitions. These definitions are from the introduction in \cite{DefantBell}.

A permutation $\pi \in S_n$ is $\mathbf{k}\textbf{-sorted}$ if and only if $\pi \in s^k(S_n)$, meaning there exists another permutation $\sigma \in S_n$ such that $s^k(\sigma) = \pi$.

For a permutation $\pi$, the function \textbf{tl}$(\pi)$ denotes the largest integer $m \in [0, n]$ such that $\pi_i = i$ for all $i$ from $n-m+1$ to $n$. For example, $13245$ has tail length $2$, while $15243$ has tail length $0$.

A \textbf{descent} of a permutation $\pi$ is an index $i$ such that $\pi_i > \pi_{i+1}$. A \textbf{descent top} is a point $(i, \pi_i)$ such that $\pi_i > \pi_{i+1}$.

A \textbf{left-to-right maximum} of a permutation $\pi$ is a point $(i, \pi_i)$ such that $\pi_j < \pi_i$ is satisfied for all indices $j < i$. Let $\textbf{LRmax}\mathbf{(\pi)}$ denote the set of left-to-right maxima of $\pi$.

$\mathbf{B_m}$ refers to the $m$th Bell number, which is the number of set partitions of $\{1, 2, \dots, m\}$ into non-empty, disjoint subsets.

$\mathbf{del^u(\pi)}$ refers to the permutation $\pi$ but with the entries $1$ through $u$ removed, then the remaining entries are standardized to $1$ to the length of $\del^u(\pi)$.

A permutation $\pi$ contains the barred pattern $\mathbf{32\overline{4}1}$ if and only if there exist indices $i_1 < i_2 < i_3$ such that both of the following hold true:
\begin{itemize}
\item $\pi_{i_1} > \pi_{i_2} > \pi_{i_3}$
\item All entries between indices $i_2$ and $i_3$ are less than $\pi_{i_1}$.
\end{itemize}
For example, $\pi = 3241$ avoids $32\overline{4}1$ because the only indices $i_1, i_2, i_3$ that satisfy the first condition are $i_1 = 1, i_2 = 2, i_3 = 4$, but this violates the second condition for containing the barred pattern $32\overline{4}1$.
\newline

\noindent
Next, here are some other results from Defant \cite{DefantBell} and Callan \cite{Callan} that we will use in our proof.
\newline

\noindent
\textbf{Theorem 2.1 \cite{Callan}:} A permutation avoids $32\overline{4}1$ if and only if every descent top of $\pi$ is a left-to-right maximum of $\pi$. Furthermore, $|\Av_n(32\overline{4}1)| = B_n$.
\newline

\noindent
\textbf{Lemma 2.2 \cite{DefantBell}:} If $\pi$ is a permutation that avoids $32\overline{4}1$ and ends in its largest entry, then there exists $\sigma \in s^{-1}(\pi)$ such that $\sigma$ avoids $32\overline{4}1$ and such that $\LRmax(\sigma) = \LRmax(\pi)$.
\newline

\noindent
\textbf{Lemma 2.3 \cite{DefantBell}:} For permutations, $s$ and $\del$ commute.
\newline

\noindent
\textbf{Lemma 2.4 \cite{DefantBell}:} Suppose $\pi \in S_n$ contains an occurrence of the pattern $32\overline{4}1$ that involves the entry $1$, and suppose $\sigma \in s^{-1}(\pi)$. Then $\sigma$ also contains an occurrence of the pattern $32\overline{4}1$ that involves the entry $1$. Moreover, there exist entries $c,d \in [3, n]$ that appear to the left of $1$ in $\sigma$ and appear to the right of $1$ in $\pi$.

For Lemma 2.4, note that when you stack sort a permutation, the entries that were on the right of $1$ remain on the right of $1$ because the entry $1$ has to leave the stack before an entry to the right of it (which is greater than $1$) can enter the stack. So stack sorting $\sigma$ increases the number of entries to the right of entry $1$ by at least $2$.
\newline

\noindent
\textbf{Lemma 2.5 \cite{BonaWalk} \cite{DefantBell} (Well Known):} Let $\pi$ be a permutation. Two entries $b, a$ form an occurrence of the pattern $21$ in $s(\pi)$ if and only if there is an entry $c$ such that $b, c, a$ form an occurrence of the pattern $231$ in $\pi$.

\section{Proof of Main Theorem}

We will first prove the if direction of the main theorem. Assume for this section that $m \geq 4$.
\newline

\noindent
\textbf{Lemma 3.1:} If a permutation $\pi \in S_{2m-4}$ satisfies $\tl(\pi) \geq m-4$ and every descent top of $\pi$ is a left-to-right maximum of $\pi$, then the permutation is $(m-4)$-sorted.
\begin{proof}
Let $\pi \in S_{2m-4}$ be a permutation that satisfies the if condition. By Theorem 2.1, because every descent top of $\pi$ is a left-to-right maximum of $\pi$, we have $\pi$ avoids $32\overline{4}1$.

Then, let $\pi'$ be $\pi$ but with the last $m-5$ entries removed. The removed entries are $m+2$ through $2m-4$, and the last entry of $\pi'$ is $m+1$, so $\pi'$ ends in its largest entry.

Since $\pi$ avoids $32\overline{4}1$, there exists no three entries $\pi_{i_1}, \pi_{i_2}, \pi_{i_3}$ in $\pi$ that form a $32\overline{4}1$ pattern. As a result, there exists no three entries $\pi_{i_1}, \pi_{i_2}, \pi_{i_3}$ in $\pi$ with $i_3 \leq m+1$ that form a $32\overline{4}1$ pattern. Since removing entries before $\pi_{i_1}$ or after $\pi_{i_3}$ doesn't affect whether $\pi_{i_1}, \pi_{i_2}, \pi_{i_3}$ forms a $32\overline{4}1$ pattern, no three entries $\pi_{i_1}, \pi_{i_2}, \pi_{i_3}$ in $\pi'$ with $i_3 \leq m+1$ form a $32\overline{4}1$ pattern. Therefore, $\pi'$ avoids $32\overline{4}1$.

Since $\pi'$ satisfies the if condition of Lemma 2.2, there exists $\sigma' \in s^{-1}(\pi')$ such that $\sigma'$ also avoids $32\overline{4}1$. Let $\sigma$ be $\sigma'$ but with the entries $m+2$ through $2m-4$ appended onto the right in increasing order. $s(\sigma') = \pi'$ implies $s(\sigma) = \pi$ because when sorting $\sigma$, all the entries $\leq m+1$ are cleared from the stack before the entry $m+2$ goes in.

$\sigma$ avoids $32\overline{4}1$ given that $\sigma'$ avoids $32\overline{4}1$ because the last $m-5$ entries are too large to be the $\pi_{i_3}$ without violating the first condition for containing $32\overline{4}1$, that $\pi_{i_1} > \pi_{i_2} > \pi_{i_3}$. We defined $\sigma$ to have $\tl(\sigma) \geq m-5$ when we appended the entries $m+2$ through $2m-4$. Therefore, we can repeat the process using $\sigma$ and $m-5$ instead of $\pi$ and $m-4$, with each repeat producing a new permutation that is the pre-image of the last with tail length reduced by at most $1$. Using this process a total of $m-4$ times eventually gives a permutation $\gamma$ such that $s^{m-4}(\gamma) = \pi$, so $\pi$ is $(m-4)$-sorted.
\end{proof}

\noindent
\textbf{Lemma 3.2:} Let $\pi \in S_{2m-4}$ be a permutation of the form $c312$, where $4 \leq c \leq m$, followed by the remaining entries placed in increasing order. Claim that $\pi$ is $(m-4)$-sorted.
\begin{proof}
Define \[\sigma = c(m+1)(m+2) \dots (2m-4)345\dots m12.\]
We create the permutation $\sigma$ by starting with $c312$, then inserting the entries $m+1$ through $2m-4$ between the $c$ and the $3$, then inserting the entries $4$ through $m$ (excluding $c$) between the $3$ and the $1$.

Now define $\sigma_i$ as $\sigma$ but with the largest $i$ entries in $4, 5, \dots, m$ (excluding $c$) and the largest $i$ entries in $(m+1), (m+2), \dots, (2m-4)$ moved to the very right. The entries in $\sigma_i$ after the entry $1$ are placed in increasing order. For example, if $m = 6$ and $\sigma = 57834612$, then $\sigma_0 = 57834612, \sigma_1 = 57341268,$ and $\sigma_2 = 53124678$. Using this definition, $\sigma = \sigma_0$ while $\pi = \sigma_{m-4}$.

We will now prove that for $0 \leq i \leq m-5$, we have $s(\sigma_i) = \sigma_{i+1}$. When we sort $\sigma_i$, first, the entries $c, m+1, m+2, \dots, 2m-5-i$ are each added into the stack and immediately popped by the next entry because the next entry is larger. Second, $2m-4-i$ is added into the empty stack. Third, the unmoved entries among $3, 4, 5, \dots, m$ excluding the largest one (let this entry be $k$) are each added into the stack and immediately popped by the next entry. Fourth, $k$ is added into the stack on top of $2m-4-i$.

The remaining entries in $\sigma_i$ form three sections: $1$ and $2$ (the first section) followed by $i$ entries that are larger than $k$ but $\leq m$ (the second section) followed by $i$ entries that are larger than $2m-4-i$ but $\leq 2m-4$ (the third section). All these remaining entries are in increasing order. In the end, these entries are outputted in increasing order, with $k$ inserted between the first and second sections and $2m-4-i$ inserted between the second and third sections. The final output is $\sigma_{i+1}$.

Therefore, $s^{m-4}(\sigma) = s^{m-4}(\sigma_0) = \sigma_{m-4} = \pi$, so $\pi$ is $(m-4)$-sorted.
\end{proof}

\noindent
In general, we have the following:

\noindent
\textbf{Lemma 3.3:} Let $\pi \in S_{2m-4}$ be a permutation of one of the forms $2c31, c231, 1c32, c132, dc21, d21$. If $c$ appears in the form, then let $4 \leq c \leq m$. If $d$ appears in the form, then let $3 \leq d \leq m$. If $c$ and $d$ both appear in the form, then let $d < c$ as well. The entries in the form are followed by the remaining entries placed in increasing order. Claim that $\pi$ is $(m-4)$-sorted.
\begin{proof}
The argument of the proof for Lemma 3.2 also works for these six forms.

Here is a table depicting the seven forms, the corresponding permutation $\sigma$, and how to create $\sigma_i$. The table includes the form from Lemma 3.2 for comparison.
\begin{center}
\begin{tabular}[t]{c | c | p{0.42\textwidth}}
Form of $\pi$ & The permutation $\sigma$ & \multicolumn{1}{c}{Creating $\sigma_i$} \\
\hline
$c312$ & $c(m+1)(m+2) \dots (2m-4)345 \dots m 12$ & \small Take $\sigma$, remove the largest $i$ entries in $(m+1)(m+2) \dots (2m-4)$ and the largest $i$ entries in $45 \dots m$ (excluding $c$), place these entries at the very end, then sort the entries after $1$ in increasing order. \\
\hline
$2c31$ & $2c(m+1)(m+2) \dots (2m-4)345 \dots m1$ & \small Take $\sigma$, remove the largest $i$ entries in $(m+1)(m+2) \dots (2m-4)$ and the largest $i$ entries in $45 \dots m$ (excluding $c$), place these entries at the very end, then sort the entries after $1$ in increasing order. \\
\hline
$c231$ & $c(m+1)(m+2) \dots (2m-4)2345 \dots m1$ & \small Take $\sigma$, remove the largest $i$ entries in $(m+1)(m+2) \dots (2m-4)$ and the largest $i$ entries in $45 \dots m$ (excluding $c$), place these entries at the very end, then sort the entries after $1$ in increasing order. \\
\hline
$1c32$ & $1c(m+1)(m+2) \dots (2m-4)345 \dots m2$ & \small Take $\sigma$, remove the largest $i$ entries in $(m+1)(m+2) \dots (2m-4)$ and the largest $i$ entries in $45 \dots m$ (excluding $c$), place these entries at the very end, then sort the entries after $2$ in increasing order. \\
\hline
$c132$ & $c(m+1)(m+2) \dots (2m-4)1345 \dots m2$ & \small Take $\sigma$, remove the largest $i$ entries in $(m+1)(m+2) \dots (2m-4)$ and the largest $i$ entries in $45 \dots m$ (excluding $c$), place these entries at the very end, then sort the entries after $2$ in increasing order. \\
\hline
$dc21$ & $dc(m+1)(m+2) \dots (2m-4)234 \dots m1$ & \small Take $\sigma$, remove the largest $i$ entries in $(m+1)(m+2) \dots (2m-4)$ and the largest $i$ entries in $34 \dots m$ (excluding $d$ and $c$), place these entries at the very end, then sort the entries after $1$ in increasing order. \\
\hline
$d21$ & \parbox[t]{0.42\textwidth}{\centering $d(m+1)(m+2) \dots (2m-4)245\dots m1e$ \\ where $e = 3$ if $d \geq 4$ and $e = 4$ if $d = 3$.} & \small Take $\sigma$, remove the largest $i$ entries in $(m+1)(m+2) \dots (2m-4)$ and the largest $i$ entries in $45 \dots m$ (excluding $d$ and $e$), place these entries at the very end, then sort the entries after $1$ in increasing order. \\
\hline
\end{tabular}

\vspace{0.25cm}
\footnotesize Table 2: Values of $\pi, \sigma, \sigma_i$ for all seven forms.
\end{center}

In each row, $\sigma_i$ is obtained by removing the largest $i$ entries from each of the two increasing blocks and moving these entries to the increasing tail at the end. When we apply $s$ to $\sigma_i$, the entries in the first increasing block (other than the largest remaining entry in the block) are pushed and immediately popped, because each such entry is followed by a larger entry. The largest remaining entry of the first increasing block remains in the stack. The same happens for the second increasing block, so after the second increasing block of $\sigma_i$ is processed, the stack contains the largest remaining entry from each block. The entries in $\sigma_i$ after the first two increasing blocks force the two entries in the stack to be outputted into the final tail so that the final tail is still in increasing order.

Since applying $s$ to $\sigma_i$ moves the largest remaining entries of the two increasing blocks to the final tail and the final tail is still increasing, $s(\sigma_i) = \sigma_{i+1}$ for $0 \leq i \leq m-5$. Therefore, $s^{m-4}(\sigma) = \pi$ because $\sigma = \sigma_0$ while $\pi = \sigma_{m-4}$.
\end{proof}

\noindent
\textbf{Lemma 3.4:} Let $\pi = 4321$. Claim that $\pi$ is $(m-4)$-sorted.
\begin{proof}
The length of the permutation is $4 = 2m-4$, so $m = 4$. It is trivial that $\pi$ is $(m-4)$-sorted because $m-4 = 0$.
\end{proof}

\noindent
We will now prove the only if direction of the main theorem. But before we begin, let us introduce another two lemmas.
\newline

\noindent
\textbf{Lemma 3.5:} Let $\pi \in S_n$. If the entries after entry $1$ in $\pi$ are in increasing order, then the entries after entry $1$ in $s(\pi)$ are also in increasing order.
\begin{proof}
Let the entries after $1$ in $\pi$ be $a_1, a_2, \dots, a_k$, where $1 < a_1 < a_2 < \dots < a_k$. While stack sorting $\pi$, let $b_1 < b_2 < \dots < b_\ell$ be the entries in the stack after all entries before entry $1$ are processed. For $1 \leq i \leq k$, let $f(i)$ be the greatest integer in $[0, \ell]$ such that $b_1, b_2, \dots, b_{f(i)}$ are all less than $a_i$.

First, the entry $1$ is entered into the stack. Next, $a_1$ pops $1, b_1, b_2, \dots, b_{f(1)}$ in this order before entering the stack. Then $a_2$ pops $a_1, b_{f(1)+1}, b_{f(1)+2}, \dots, b_{f(2)}$ in this order before entering the stack. This continues until $a_k$ pops $a_{k-1}, b_{f(k-1)+1}, b_{f(k-1)+2}, \dots, b_{f(k)}$ in this order before entering the stack. The output so far (starting with the entry $1$) is $1, a_1, a_2, \dots, a_{k-1}, b_1, b_2, \dots, b_{f(k)}$ in increasing order. Now that there are no entries left in $\pi$, the remaining entries in the stack, $a_k, b_{f(k)+1}, b_{f(k)+2}, \dots, b_\ell,$ are outputted in increasing order. Since $a_k > b_{f(k)}$, the entries in $s(\pi)$ after the entry $1$ are indeed in increasing order.
\end{proof}

\noindent
\textbf{Lemma 3.6:} Let $\pi \in S_n$ be $k$-sorted, where $0 \leq k \leq n$. Then $\tl(\pi) \geq k$.
\begin{proof}
We will prove this lemma using induction.

Base Case $k=0$: $\tl(\pi)$ is nonnegative by definition, so the lemma is true for $k = 0$.

Inductive Hypothesis: For some $0 \leq k \leq n-1$, if a permutation $\pi \in S_n$ is $k$-sorted, then $\tl(\pi) \geq k$.

Inductive Step: Take any permutation $\pi \in S_n$ that is $(k+1)$-sorted. Then there exists $\sigma \in S_n$ such that $\sigma$ is $k$-sorted and $s(\sigma) = \pi$. By the inductive hypothesis, $\tl(\sigma) \geq k$. When we stack sort $\sigma$, the entry $n-k$ sits at the bottom of the stack until the entry $n-k+1$ pops it. Because the entries $n-k+1$ through $n$ are in increasing order in $\sigma$, each one pops the one immediately before, so $n-k$ through $n$ appear in increasing order at the end of $\pi$. Since we proved that $\tl(\pi) \geq k+1$ for any permutation $\pi \in S_n$ that is $(k+1)$-sorted, the inductive step is complete.
\end{proof}

\noindent
\textbf{Claim 3.7:} Let $\pi \in S_{2m-4}$ be $(m-4)$-sorted and contain a $32\overline{4}1$ pattern. Then, $\pi$ satisfies one of the following:
\begin{itemize}
    \item The first four (three for the last case) entries of $\pi$ are of the form $c312, 2c31, c231, 1c32, c132, dc21,$ or $d21$, where $4 \leq c \leq m$ and $3 \leq d \leq m$ and $d < c$ (each of the three is used if applicable), and the remaining entries are placed in increasing order.
    \item $\pi = 4321$.
\end{itemize}
\begin{proof}
The fact that $\pi$ is $(m-4)$-sorted means there exists $\sigma \in S_{2m-4}$ such that $\pi = s^{m-4}(\sigma)$. By Lemma 2.3, we get
\[\del(\pi) = \del(s^{m-4}(\sigma)) = s^{m-4}(\del(\sigma)).\]
Since $m-4 = (m-1)-3$ and the length of $\del(\pi)$ and $\del(\sigma)$ is $2m-5 = 2(m-1)-3$, we use Theorem 1.2 to get that either:
\begin{itemize}
\item $\tl(\del(\pi)) \geq m-4$ and every descent top of $\del(\pi)$ is a left-to-right maximum of $\del(\pi)$.
\item Or $\del(\pi) = \zeta_{\ell, m-1}$ for some $\ell \in \{3, \dots, m-1\}$. The permutation $\zeta_{\ell, m-1}$ is obtained from the identity permutation of length $2m-5$ by swapping entries $1$ and $2$ and moving entry $\ell$ to the front.
\end{itemize}
So we can split this claim.
\newline

\noindent
\textbf{Subclaim 3.7.1:} Let $\pi \in S_{2m-4}$ be $(m-4)$-sorted and contain a $32\overline{4}1$ pattern. Also suppose $\tl(\del(\pi)) \geq m-4$ and every descent top of $\del(\pi)$ is a left-to-right maximum of $\del(\pi)$. Then, $\pi$ satisfies one of the following:
\begin{itemize}
    \item The first four (three for the last case) entries of $\pi$ are of the form $c312, 2c31, c231, 1c32, c132, dc21,$ or $d21$, where $4 \leq c \leq m$ and $3 \leq d \leq m$ and $d < c$ (each of the three is used if applicable), and the remaining entries are placed in increasing order.
    \item $\pi = 4321$.
\end{itemize}
\begin{subproof}
By Theorem 2.1, not every descent top of $\pi$ is a left-to-right maximum of $\pi$. However, we are given that every descent top of $\del(\pi)$ is a left-to-right maximum of $\del(\pi)$. Since $\pi$ is obtained from $\del(\pi)$ by increasing all entries by $1$ and then inserting the entry $1$, the only possible new descent top of $\pi$ is the entry immediately before $1$. Therefore, the descent top of $\pi$ that is not a left-to-right maximum of $\pi$ must be the entry immediately before $1$.
\begin{itemize}
\item If the entry $1$ is placed as the first entry, then the set of descent tops and the set of left-to-right maxima are unchanged other than the fact that $1$ is added to the latter set.
\item If the entry $1$ is placed after a descent top, then neither set is affected.
\item If the entry $1$ is placed after an entry that is not a descent top but is a left-to-right maximum, then even though the entry before $1$ becomes a descent top, all descent tops are still left-to-right maxima.
\item If the entry $1$ is placed after an entry that is not a descent top and is not a left-to-right maximum, the entry becomes a descent top but still is not a left-to-right maximum. This is the only case in which not every descent top of $\pi$ is a left-to-right maximum of $\pi$.
\end{itemize}
Therefore, $\pi$ must contain a $32\overline{4}1$ pattern involving the entry $1$: let $\pi_{i_3} = 1,$ let $\pi_{i_2}$ be the entry immediately before, and let $\pi_{i_1}$ be the entry to the left of $\pi_{i_2}$ that prevents it from being a left-to-right maximum.

Because $\pi$ is $(m-4)$-sorted and contains a $32\overline{4}1$ pattern involving the entry $1$, there exist permutations \[\pi_1, \pi_2, \dots, \pi_{m-4}\] such that $s(\pi_i) = \pi_{i-1}$ for $1 \leq i \leq m-4$ and $\pi = \pi_0$. Let $\sigma = \pi_{m-4}$.

Since $\pi$ contains a $32\overline{4}1$ pattern involving the entry $1$, we can use Lemma 2.4 to get that $\pi_1$ also contains a $32\overline{4}1$ pattern involving the entry $1$. Then apply it repeatedly to get that $\pi_2, \pi_3, \dots, \pi_{m-4}$ also each contain a $32\overline{4}1$ pattern involving the entry $1$.

Each time Lemma 2.4 is applied, at least two entries move from the left of entry $1$ to the right of entry $1$ when $s$ is used on $\pi_i$. However, no entries move from the right of entry $1$ to the left. If there exists an entry to the right of entry $1$, then entry $1$ is immediately popped (by the next entry, which is always larger) before the entries after it can enter the stack.

Therefore, after $m-4$ applications of Lemma 2.4, we get that $\pi$ has at least $2m-8$ entries $\geq 3$ to the right of entry $1$. But since the entry $1$ is the $\pi_{i_3}$ of the $32\overline{4}1$ pattern, there exists an entry $\geq 3$ to the left of entry $1$, so there can be at most $2m-7$ entries $\geq 3$ to the right of entry $1$ in $\pi$.

Because there are at least two entries to the left of entry $1$ in $\pi$, there are at least $2m-6$ entries to the left of entry $1$ in $\sigma = \pi_{m-4}$, so there can be at most $1$ entry to the right of entry $1$ in $\sigma$. The entries after entry $1$ in $\sigma$ are in increasing order, and by Lemma 3.5, the entries after entry $1$ in $\pi$ are also in increasing order.

Also, there is either $1$ or $2$ entries $\geq 3$ to the left of entry $1$ in $\pi$, since there are $2m-6$ entries $\geq 3$ in $\pi$ in total. Therefore, we have the following three sub-cases:
\begin{itemize}
\item If there is only one entry $\geq 3$ to the left of entry $1$ in $\pi$, then the entry $2$ must serve as the $\pi_{i_2}$ in the $32\overline{4}1$ pattern involving the entry $1$. Therefore, the only form for $\pi$ allowed in this sub-case is indeed the $d21$ form mentioned in the first bullet point of this claim.
\item If there are two entries $\geq 4$ to the left of entry $1$ in $\pi$, let them be $c, d,$ where $d < c$. Assume for the sake of contradiction that $c$ appears before $d$. The entries $c, d, 1$ form a $32\overline{4}1$ pattern (no entries larger than $c$ are on the left of $1$, let alone between $d$ and $1$).

Let $\sigma^{(1)} \in s^{-1}(\pi)$. By Lemma 2.5, there must be entries $c', d'$ in $\sigma^{(1)}$ such that $c' > c, d' > d$ and $c'$ is in between $c$ and $d$ while $d'$ is in between $d$ and $1$. The entries on the left of $1$ in $\pi$ must still be on the left of $1$ in $\sigma^{(1)}$ (see the note for Lemma 2.4), so at least two more entries with value $\geq 4$ are to the left of entry $1$ in $\sigma^{(1)}$ when compared to $\pi$.

Apply this process $m-5$ more times to get that for $\sigma^{(m-4)} \in s^{-(m-4)}(\pi)$, there exist at least $2m-8$ more entries with value $\geq 4$ that are to the left of $1$ in $\sigma^{(m-4)}$ when compared to $\pi$. This is a contradiction, since $\pi$ is a standardized length $2m-4$ permutation, meaning it does not have $2m-6$ entries with value $\geq 4$.

So $d, c, 1$ appear in $\pi$ in this order. In order for $\pi$ to contain a $32\overline{4}1$ pattern involving the entry $1$, the entry $2$ must appear specifically between $c$ and $1$. Therefore, the only form for $\pi$ allowed in this sub-case is indeed the $dc21$ form mentioned in the first bullet point of this claim.
\item The remaining sub-case is: there is one entry $\geq 4$ and the entry $3$ to the left of entry $1$ in $\pi$. Let $c$ be the entry $\geq 4$ that is before $1$.

If $c$ is placed after $3$ in $\pi$, then in order for $\pi$ to contain a $32\overline{4}1$ pattern involving the entry $1$, the entry $2$ must appear specifically between $c$ and $1$, giving us the $3c21$ form mentioned in the first bullet point of this claim.

Otherwise, $c$ is placed before $3$ in $\pi$. In order for the entry $2$ to be in between $3$ and $1$ in $\pi$, by Lemma 2.5, for any $\sigma^{(1)} \in s^{-1}(\pi)$, there exists $a', b', c'$ such that $c' > c, b' > 3, a' > 2$ and $c'$ is between $c$ and $3$ while $b'$ is between $3$ and $2$ while $a'$ is between $2$ and $1$. Apply this process $m-5$ more times to get that for $\sigma^{(m-4)} \in s^{-(m-4)}(\pi)$, there exist at least $3(m-4)+2$ entries $\geq 3$ to the left of entry $1$ in $\sigma^{(m-4)}$, which in total has $2m-6$ entries $\geq 3$. So we would need \[2m-6 \geq 3(m-4)+2 = 3m-10 \implies m \leq 4\]
Hence, if $c, 3, 2, 1$ are placed in $\pi$ in this order, then we must have $m = 4$ and $c = 4$. So the only forms for $\pi$ allowed in this sub-sub-case, where $c$ is before $3$ in $\pi$, are the $c312, 4321, c231, 2c31$ forms mentioned in the first and second bullet points of this claim.
\end{itemize}
\end{subproof}

\noindent
\textbf{Claim 3.7.2:} Let $\pi \in S_{2m-4}$ be $(m-4)$-sorted and contain a $32\overline{4}1$ pattern. Also let $\del(\pi) = \zeta_{\ell, m-1}$ for some $\ell \in \{3, \dots, m-1\}$. The permutation $\zeta_{\ell, m-1}$ is obtained from the identity permutation of length $2m-5$ by swapping entries $1$ and $2$ and moving entry $\ell$ to the front. Then, $\pi$ satisfies one of the following:
\begin{itemize}
    \item The first four (three for the last case) entries of $\pi$ are of the form $c312, 2c31, c231, 1c32, c132, dc21,$ or $d21$, where $4 \leq c \leq m$ and $3 \leq d \leq m$ and $d < c$ (each of the three is used if applicable), and the remaining entries are placed in increasing order.
    \item $\pi = 4321$.
\end{itemize}
\begin{subproof}
If the \textbf{second bullet point} about $\del(\pi)$ is true, then $\pi$ is of the form $c3245 \dots (2m-4)$ with the entry $1$ added somewhere in here.
\begin{itemize}
\item If we place $1$ before $c$, then we get the $1c32$ form.
\item If we place $1$ between $c$ and $3$, then we get the $c132$ form.
\item If we place $1$ between $3$ and $2$, then we get the $c312$ form.
\item If we place $1$ right after $2$, we get the same $c,3,2,1$ scenario which we proved a few paragraphs above.
\item If we place $1$ after $2$ such that there are $k \geq 1$ entries in between them, then the proof is very similar to the previous case because there are still three descents to the left of entry $1$ in $\pi$ (the descents are $c,3,$ and the entry immediately before $1$). However, there are initially $k+2$ (instead of two) entries $\geq 3$ to the left of entry $1$ in $\pi$, so the inequality becomes \[2m-6 \geq 3(m-4)+(k+2) = 3m-10+k \implies m \leq 4-k.\]
This is a contradiction because we assumed $m \geq 4$, so we get no form for this sub-case.
\end{itemize}
\end{subproof}
We proved that as long as $\pi \in S_{2m-4}$ is $(m-4)$-sorted and contains a $32\overline{4}1$ pattern, then $\pi$ satisfies one of the following no matter what:
\begin{itemize}
    \item The first four (three for the last case) entries of $\pi$ are of the form $c312, 2c31, c231, 1c32, c132, dc21,$ or $d21$, where $4 \leq c$ and $3 \leq d$ and $d < c$ (each of the three is used if applicable), and the remaining entries are placed in increasing order.

    $c,d \leq m$ comes from Lemma 3.6, which tells us that the entries $m+1$ through $2m-4$ in $\pi$ are listed in increasing order at the end of $\pi$ because $\pi$ is $(m-4)$-sorted.
    \item $\pi = 4321$.
\end{itemize}
Therefore, our proof of Claim 3.7 is complete.
\end{proof}
We will now prove the enumeration part of the main theorem.
\newline

\noindent
\textbf{Claim 3.8:} If $m \geq 5,$ then
\[|s^{m-4}(S_{2m-4})| = B_m + \frac{m^2 + 7m - 28}{2}.\]
Add one to the above value if $m = 4$.
\begin{proof}
For $m = 4$, we have $s^{m-4}(S_{2m-4}) = s^0(S_4) = S_4$, so the left side is equal to $24$. The right side is equal to $15 + 8 + 1 = 24$, so the equation is true for $m = 4$.

For $m > 4$, we will first prove that the three bullet points of Theorem 1.3 do not overlap.

The only permutation in the third bullet point does not satisfy $m > 4$.

The first and second bullet points do not overlap. In the $c312, 2c31, c231, 1c32,$ and $c132$ forms, the $3$ is a descent top but not a left-to-right maximum. In the $dc21$ and $d21$ forms, the $2$ is a descent top but not a left-to-right maximum. So any permutation that satisfies the second bullet point does not satisfy the first bullet point.

For the first bullet point, we will construct a bijection between:
\begin{itemize}
\item The permutations $\pi \in S_{2m-4}$ such that $\tl(\pi) \geq m-4$ and every descent top of $\pi$ is a left-to-right maximum of $\pi$.
\item The set partitions of $1, 2, \dots, m$.
\end{itemize}
To get from the first to the second, construct $\pi'$ by removing the last $m-4$ entries in $\pi$. Then, take each left-to-right maximum in $\pi'$ and combine it with the entries between it and the next left-to-right maximum to form a set. These sets form a set partition of $1, 2, \dots, m$.

This process is reversible. Given a set partition of $\{1, 2, \dots, m\}$, first arrange the sets such that in each set, the largest element appears at the front and the rest are arranged in increasing order. For example, $\{2, 5, 1, 3\}$ turns into $\{5, 1, 2, 3\}$. Then, sort the sets in increasing order by their first elements and merge the sets in this order. Finally, add $m+1$ through $2m-4$ in increasing order on the right to get the permutation. The descent tops of the resulting permutation are the first elements of the sorted sets with $\geq 2$ elements, while the left-to-right maxima of the resulting permutation are the first elements of all of the sorted sets, so the descent-top condition is satisfied.

For example, take the permutation $\pi = \{5, 1, 2, 3, 7, 4, 6, 8, 9, 10\}$ with $m = 7$. Then $\pi' = \{5, 1, 2, 3, 7, 4, 6\}$, which has left-to-right maxima $5$ and $7$, so the sets are $\{5, 1, 2, 3\}$ and $\{7, 4, 6\}$. Now, if we take these sets and sort them, $\{5, 1, 2, 3\}$ comes before $\{7, 4, 6\}$ because $5 < 7$. Merging the sets gives $\{5, 1, 2, 3, 7, 4, 6\}$. And adding $m+1 = 8$ through $2m-4 = 10$ gives $\pi = \{5, 1, 2, 3, 7, 4, 6, 8, 9, 10\}$ again.

Therefore, there are $B_m$ permutations $\pi \in S_{2m-4}$ such that $\tl(\pi) \geq m-4$ and every descent top of $\pi$ is a left-to-right maximum of $\pi$.

For the second bullet point, we can easily tell that the forms are mutually exclusive. In the first five forms, there are $m-3$ ways to select $c$, so there are $m-3$ permutations of each form. In the sixth form, there are $\binom{m-2}{2}$ ways to select $c$ and $d$. In the seventh form, there are $m-2$ ways to select $d$. Therefore, there are \[5(m-3)+\binom{m-2}{2} + (m-2) = \frac{m^2 + 7m - 28}{2}\]
permutations satisfying the second bullet point of Theorem 1.3.

And there are $0$ permutations satisfying the third bullet point of Theorem 1.3.

Adding all together,\[|s^{m-4}(S_{2m-4})| = B_m + \frac{m^2 + 7m - 28}{2}\] for $m > 4$.
\end{proof}

\section{General Remarks and Future Work}

Defant posed the following question to the author (privately):

\noindent
\textbf{Question 4.1:} In general, for fixed $k \geq 2$, is there a characterization of the following form?

If $\pi \in S_{2m-k}$ is $(m-k)$-sorted, then one of the following holds:
\begin{itemize}
\item $\tl(\pi) \geq n-m$, where $n$ is the length of $\pi$, and every descent top of $\pi$ is a left-to-right maximum of $\pi$.
\item The entries after the first $k$ entries are all in increasing order.
\end{itemize}

\noindent
\textit{Proof:} The answer is no. There is a characterization for $k = 2, 3, 4$, but not a characterization for $k = 5$.

This is true for $k = 2$ because all permutations $\pi \in S_{2m-2}$ that are $(m-2)$-sorted satisfy the first bullet point. This is true for $k = 3$ because the permutations $\pi \in S_{2m-3}$ that are $(m-3)$-sorted either satisfy the first bullet point or are of the $d21$ form mentioned in Theorem 1.3. This is also true for $k = 4$ because the forms $c312, 2c31, c231, 1c32, c132, dc21, d21,$ and $4321$ all have the entries after the first four entries in increasing order.

However, this is false for $k = 5$. A counterexample found using computer code is $\pi = 6, 2, 1, 3, 4, 7, 5, 8, 9, 10, 11$. We have $m = 8$ and $\sigma = 6, 7, 8, 11, 2, 3, 4, 9, 1, 10, 5$, where $s^{m-5}(\sigma) = \pi$.

Since there exists $\sigma \in S_{2m-5}$ such that $s^{m-5}(\sigma) = \pi$, the permutation $\pi$ is $(m-5)$-sorted. The entry $2$ is a descent top of $\pi$ but not a left-to-right maximum of $\pi$, so the first bullet point is not satisfied. The entries $7, 5, 8, 9, 10, 11$ in $\pi$ are not in increasing order, so the second bullet point is not satisfied either.
\hfill $\square$ \newline

An implication of the answer to this question is that the proof for the $2m-5$ case would be much more difficult and lengthy compared to the $2m-4$ case, since there can be forms of length $> 5$.

This is not the only important property that holds for $2m-3$ and $2m-4$ but not $2m-5$. In Defant's proof of the $2m-3$ case and our proof of the $2m-4$ case, we used this property: the entries after entry $1$ in $\pi \in S_{2m-k}$ are in increasing order, given that $\pi$ is $(m-k)$-sorted and $\tl(\del(\pi)) \geq m-k$ and every descent top of $\del(\pi)$ is a left-to-right maximum of $\del(\pi)$.

For the $2m-5$ case, here is a counterexample: \[\pi = 4, 3, 1, 5, 6, 2, 7, 8, 9, 10, 11.\]
The $\pi$ in this counterexample satisfies all three givens:
\begin{itemize}
\item $\pi$ is $(m-5)$-sorted, because $\sigma = 4, 6, 9, 11, 3, 5, 7, 8, 1, 10, 2$ satisfies $\sigma \in S_{11}$ and $s^3(\sigma) = \pi$.
\item We have $\del(\pi) = 3, 2, 4, 5, 1, 6, 7, 8, 9, 10$. The value $\tl(\del(\pi)) = 5 \geq 3 = m-k$.
\item Every descent top of $\del(\pi)$ (which are the entries $3$ and $5$) is a left-to-right maximum of $\del(\pi)$ (which are the entries $\geq 3$).
\end{itemize}
But the entries after $1$ in $\pi \in S_{2m-5}$ are not in increasing order (the $2$ should come before the $5$).

Therefore, work on the $2m-5$ case would help reveal the similarities and differences between the behaviors of the $2m-3$ and $2m-4$ cases compared to the behaviors of the $2m-5$ and later cases.
\newline

So far, for $k = 3$ and $k = 4$, the value $|s^{m-k}(S_{2m-k})|$ is equal to $B_m$ plus a degree $k-2$ polynomial. It is desirable to understand whether this pattern continues, so we pose the following conjecture:
\newline

\noindent
\textbf{Conjecture 4.2:} For each positive integer $k \geq 3$, there exists a finite positive integer $c$ and a degree $k-2$ polynomial $p(m)$ such that for all $m \geq c,$ \[|s^{m-k}(S_{2m-k})| = B_m + p(m).\]

\noindent
If this conjecture is true, it would also be desirable to understand the growth of the least possible $c$ as $k$ increases. For instance, the least possible $c$'s for $k = 3$ and $k = 4$ are $3$ and $5$, respectively.

\begin{center}
\large
$\textbf{Acknowledgments}$
\end{center}

The author gratefully acknowledges Colin Defant for introducing the author to the open questions in his 2022 thesis \textit{Stack-Sorting and Beyond} \cite{DefantThesis}, including the $2m-4$ case of highly sorted permutations. The author also thanks Defant for identifying places in earlier drafts where the rigor and clarity could be improved. The author developed the results and proofs in this paper.


\begin{thebibliography}{10}

\bibitem{Knuth}
Knuth, Donald E. \textit{The Art of Computer Programming}. 1968. Vol. 1, Addison-Wesley Professional, 1997, pp. 238–243.

\bibitem{West}
West, Julian. \textit{Permutations with Forbidden Subsequences, And, Stack-Sortable Permutations}. 1990.

\bibitem{BonaWalk}
B\'ona, Mikl\'os. \textit{A Walk through Combinatorics}. 5th ed., World Scientific, 2023, pp. 361–371.

\bibitem{Zeilberger}
Doron Zeilberger. “A Proof of Julian West’s Conjecture That the Number of Two-Stacksortable Permutations of Length N Is 2(3n)!/((N + 1)!(2n + 1)!).” \textit{Discrete Mathematics}, vol. 102, no. 1, May 1992, pp. 85–93, https://doi.org/10.1016/0012-365X(92)90351-F.

\bibitem{Bousquet}
Bousquet-Mélou, Mireille. “Sorted And/or Sortable Permutations.” \textit{Discrete Mathematics}, vol. 225, no. 1-3, 2000, pp. 25–50, https://doi.org/10.1016/S0012-365X(00)00146-1.

\bibitem{Claesson}
Claesson, Anders, and Henning Úlfarsson. “Sorting and Preimages of Pattern Classes.” \textit{Discrete Mathematics \& Theoretical Computer Science}, vol. DMTCS Proceedings vol. AR, 24th International Conference on Formal Power Series and Algebraic Combinatorics (FPSAC 2012), 2012, https://doi.org/10.46298/dmtcs.3066.

\bibitem{DefantDescents}
Defant, Colin. “Descents in $t$-Sorted Permutations.” \textit{Journal of Combinatorics}, vol. 11, no. 3, 2020, pp. 511–526, https://doi.org/10.4310/joc.2020.v11.n3.a5.

\bibitem{DefantBell}
Defant, Colin. “Highly Sorted Permutations and Bell Numbers.” \textit{Enumerative Combinatorics and Applications}, vol. 2021, no. 1, 4 Dec. 2020, https://doi.org/10.54550/ECA2021V1S1R6.

\bibitem{Callan}
Callan, David. “A Combinatorial Interpretation of the Eigensequence for Composition.” \textit{Journal of Integer Sequences}, vol. 9, no. 1, 2006.

\bibitem{DefantThesis}
Defant, Colin. \textit{Stack-Sorting and Beyond}. 2022, arks.princeton.edu/ark:/88435/dsp016m311s469.

\end{thebibliography}
\end{document}